\documentclass[11pt]{amsart}

\usepackage[T1]{fontenc}
\usepackage[utf8]{inputenc}
\usepackage{lmodern}
\usepackage{microtype}
\usepackage{amsmath,amssymb,amsthm,mathtools}
\usepackage{tikz}
\usepackage{tikz-cd}
\usepackage{enumitem}
\usepackage{hyperref}
\usepackage[nameinlink,capitalise,noabbrev]{cleveref}

\hypersetup{
  colorlinks=true,
  linkcolor=blue,
  citecolor=blue,
  urlcolor=blue
}

\newtheorem{theorem}{Theorem}[section]
\newtheorem{lemma}[theorem]{Lemma}
\newtheorem{corollary}[theorem]{Corollary}

\theoremstyle{definition}
\newtheorem{definition}[theorem]{Definition}

\theoremstyle{remark}
\newtheorem{remark}[theorem]{Remark}

\crefname{theorem}{Theorem}{Theorems}
\Crefname{theorem}{Theorem}{Theorems}
\crefname{lemma}{Lemma}{Lemmas}
\Crefname{lemma}{Lemma}{Lemmas}
\crefname{corollary}{Corollary}{Corollaries}
\Crefname{corollary}{Corollary}{Corollaries}
\crefname{proposition}{Proposition}{Propositions}
\Crefname{proposition}{Proposition}{Propositions}
\crefname{definition}{Definition}{Definitions}
\Crefname{definition}{Definition}{Definitions}
\crefname{example}{Example}{Examples}
\Crefname{example}{Example}{Examples}
\crefname{remark}{Remark}{Remarks}
\Crefname{remark}{Remark}{Remarks}

\newcommand{\kk}{\Bbbk}
\newcommand{\A}{\mathbb A}

\newcommand{\Ocal}{\mathcal O}
\newcommand{\Rep}{\operatorname{Rep}}
\newcommand{\rep}{\operatorname{rep}}
\newcommand{\End}{\operatorname{End}}
\newcommand{\Hom}{\operatorname{Hom}}
\newcommand{\Ext}{\operatorname{Ext}}
\newcommand{\GL}{\operatorname{GL}}

\newcommand{\depth}{\operatorname{depth}}

\newcommand{\pd}{\operatorname{pd}}

\title[Quiver Orbit Closures]{A Homological Approach to the Cohen--Macaulayness Problem for Quiver Orbit Closures}
\author{Ales Bouhada}
\address{Department of Mathematics, University of Sherbrooke, Sherbrooke, Quebec, Canada}
\email{alesm.bouhada@gmail.com}
\subjclass[2010]{16E10, 16E65, 16G20, 14M12, 13C14}
\keywords{Quiver varieties, orbit closures, Cohen--Macaulay varieties, projective dimension, representations of quivers}

\begin{document}

\begin{abstract}
We propose a homological approach to the Cohen--Macaulayness problem for orbit closures in varieties of quiver representations.  Given a finite quiver $Q$, a dimension vector $d$, and a representation $M\in \rep(Q,d)$, we relate the Cohen--Macaulay property of the orbit closure $\overline{\Ocal_M}$ to the projective dimension of its coordinate ring as a module over the polynomial algebra $\kk[\rep(Q,d)]$.  For tree quivers, we obtain a numerical criterion expressing Cohen--Macaulayness in terms of the projective dimension of $\kk[\overline{\Ocal_M}]$, the dimension of the representation variety, and the dimension of the endomorphism algebra of $M$.  This viewpoint recasts known results for Dynkin quivers of types $A_n$ and $D_n$ and suggests a possible homological strategy for the remaining exceptional Dynkin types.
\end{abstract}

\maketitle

\section{Introduction}

Let $\kk$ be an algebraically closed field.  The geometry of orbit closures in varieties of quiver representations is a central theme in the representation theory of finite-dimensional algebras.  Let $Q=(Q_0,Q_1,s,t)$ be a finite quiver and let $d=(d_i)_{i\in Q_0}$ be a dimension vector.  The affine variety
\[
\rep(Q,d)=\prod_{\alpha\in Q_1}\Hom_\kk(\kk^{d_{s(\alpha)}},\kk^{d_{t(\alpha)}})
\]
parametrizes representations of $Q$ with dimension vector $d$.  The algebraic group
\[
G_d=\prod_{i\in Q_0}\GL_{d_i}(\kk)
\]
acts on $\rep(Q,d)$ by change of bases.  Its orbits are precisely the isomorphism classes of representations with dimension vector $d$.

For a representation $M\in \rep(Q,d)$, we denote by $\Ocal_M$ its orbit and by $\overline{\Ocal_M}$ its Zariski closure.  The geometry of $\overline{\Ocal_M}$ reflects subtle representation-theoretic information about $M$ and its degenerations.  In particular, one is interested in whether these orbit closures are normal, Cohen--Macaulay, Gorenstein, or have rational singularities.

For quivers of type $A_n$ and $D_n$, orbit closures are known to be normal and Cohen--Macaulay with rational singularities; see \cite{AbeasisDelFraKraft1981,BobinskiZwara2001,BobinskiZwara2002}.  For the exceptional Dynkin quivers of types $E_6$, $E_7$, and $E_8$, the corresponding Cohen--Macaulayness problem is more delicate.  The purpose of this paper is to explain how one may attack this problem using homological algebra.

Let
\[
R=\kk[\rep(Q,d)]
\]
be the coordinate ring of the representation variety.  Since $R$ is a polynomial algebra, Hilbert's syzygy theorem implies that the $R$-module $\kk[\overline{\Ocal_M}]$ has finite projective dimension.  The main point of this paper is that, for tree quivers, the Cohen--Macaulay property of $\overline{\Ocal_M}$ is equivalent to an explicit formula for this projective dimension.  The formula involves the dimension of the representation variety and the endomorphism algebra of $M$.

The resulting criterion is particularly useful because it transforms a geometric question into a homological one.  In favorable situations, the projective dimension of $\kk[\overline{\Ocal_M}]$ can be approached through a minimal free resolution of the defining ideal of the orbit closure.  Thus the Cohen--Macaulayness problem becomes closely related to the syzygies of orbit closures.

\section{Preliminaries}

\subsection{Quiver representations and representation varieties}

A quiver is a quadruple
\[
Q=(Q_0,Q_1,s,t),
\]
where $Q_0$ is the set of vertices, $Q_1$ is the set of arrows, and
\[
s,t\colon Q_1\longrightarrow Q_0
\]
assign to each arrow its source and target.  A quiver is called a tree quiver if its underlying graph is connected and has no cycles.  In particular, every Dynkin quiver is a tree quiver.

A representation $M$ of $Q$ consists of a family of finite-dimensional $\kk$-vector spaces $(M_i)_{i\in Q_0}$ together with linear maps
\[
M_\alpha\colon M_{s(\alpha)}\longrightarrow M_{t(\alpha)},
\qquad \alpha\in Q_1.
\]
A morphism $h\colon M\to N$ of representations is a family of linear maps
\[
h_i\colon M_i\longrightarrow N_i,
\qquad i\in Q_0,
\]
such that, for every arrow $\alpha\colon i\to j$, the diagram
\[
\begin{tikzcd}
M_i \arrow[r,"M_\alpha"] \arrow[d,"h_i"'] & M_j \arrow[d,"h_j"] \\
N_i \arrow[r,"N_\alpha"'] & N_j
\end{tikzcd}
\]
commutes.  We denote by $\Hom_Q(M,N)$ the vector space of morphisms from $M$ to $N$, and by $\Rep_\kk(Q)$ the category of finite-dimensional representations of $Q$.  This category is equivalent to the category of finite-dimensional modules over the path algebra $\kk Q$.

Fix a dimension vector $d=(d_i)_{i\in Q_0}$.  Then
\[
\rep(Q,d)=\prod_{\alpha\in Q_1}\Hom_\kk(\kk^{d_{s(\alpha)}},\kk^{d_{t(\alpha)}})
\]
is an affine space of dimension
\[
\ell=\sum_{\alpha\in Q_1}d_{s(\alpha)}d_{t(\alpha)}.
\]
We shall write $\rep(Q,d)\cong \A^\ell$.

The group
\[
G_d=\prod_{i\in Q_0}\GL_{d_i}(\kk)
\]
acts on $\rep(Q,d)$ by
\[
(g_i)_{i\in Q_0}\cdot (M_\alpha)_{\alpha\in Q_1}
=
(g_{t(\alpha)}M_\alpha g_{s(\alpha)}^{-1})_{\alpha\in Q_1}.
\]
The orbit $\Ocal_M$ of a point $M\in \rep(Q,d)$ consists precisely of the representations isomorphic to $M$.

The coordinate ring of $\rep(Q,d)$ is the polynomial algebra
\[
R=\kk[\rep(Q,d)]
 =
\kk[x_{ij}^{\alpha}\mid \alpha\in Q_1,
1\leq i\leq d_{t(\alpha)},\ 1\leq j\leq d_{s(\alpha)}].
\]
If $I(\overline{\Ocal_M})$ denotes the ideal of functions vanishing on $\overline{\Ocal_M}$, then
\[
\kk[\overline{\Ocal_M}]
\cong
R/I(\overline{\Ocal_M}).
\]

\subsection{Degenerations and extension spaces}

Let $M,N\in \rep(Q,d)$.  We say that $M$ degenerates to $N$, and write
\[
M\leq_{\deg}N,
\]
if $N\in \overline{\Ocal_M}$.  Thus degenerations are precisely the points lying in orbit closures.

We recall two standard facts connecting the geometry of orbits with representation theory.  First,
\[
\Ocal_M \text{ is open in } \rep(Q,d)
\quad\Longleftrightarrow\quad
\Ext_Q^1(M,M)=0.
\]
Second,
\[
\Ocal_M \text{ is closed}
\quad\Longleftrightarrow\quad
M \text{ is semisimple}.
\]

We shall also use the Artin--Voigt formula.

\begin{theorem}[Artin--Voigt]
Let $M\in \rep(Q,d)$. Then
\[
\dim \rep(Q,d)-\dim \Ocal_M
=
\dim_\kk \Ext_Q^1(M,M).
\]
\end{theorem}

\subsection{Commutative algebra background}

Let
\[
S=\kk[x_1,\dots,x_\ell]
\]
be a polynomial algebra, let \(I\subseteq S\) be a homogeneous ideal, and set
\[
R=S/I.
\]
Let \(L\) be a finitely generated graded \(R\)-module. We denote by
\(\pd_S L\) its projective dimension as an \(S\)-module. We say that \(L\) is
Cohen--Macaulay if
\[
\depth_{\mathfrak m}L=\dim L,
\]
where
\[
\mathfrak m=(x_1,\dots,x_\ell)/I
\]
is the irrelevant maximal ideal of \(R\). An affine variety \(X\) is
Cohen--Macaulay if its coordinate ring \(\kk[X]\) is Cohen--Macaulay.

We shall use the Auslander--Buchsbaum formula.

\begin{theorem}[Auslander--Buchsbaum]
Let
\[
S=\kk[x_1,\dots,x_\ell]
\]
be a polynomial algebra, and let \(L\) be a finitely generated graded
\(S\)-module. Then
\[
\pd_S L+\depth_{(x_1,\dots,x_\ell)}L=\dim S.
\]
\end{theorem}

\section{Main results}

\begin{definition}
Let $M\in \rep(Q,d)$.  We say that $M$ is \emph{homogeneous} if the orbit closure $\overline{\Ocal_M}$ is an affine cone, equivalently, if it contains every line spanned by one of its points.
\end{definition}

For $\lambda\in \kk$, we denote by $\lambda M$ the representation corresponding to the point $\lambda m\in \rep(Q,d)$, where $m$ is the point representing $M$.

\begin{lemma}\label{lem:homogeneous}
A representation $M$ is homogeneous if and only if
\[
\lambda M\cong M
\]
for every $\lambda\in \kk^\times$.
\end{lemma}

\begin{proof}
Assume first that $M$ is homogeneous.  Then $\lambda m\in \overline{\Ocal_M}$ for every $\lambda\in\kk$.  Since orbit closures are noetherian and orbits are constructible, the descending chain
\[
\overline{\Ocal_M}
\supseteq
\overline{\Ocal_{\lambda M}}
\supseteq
\overline{\Ocal_{\lambda^2M}}
\supseteq \cdots
\]
stabilizes.  Hence, for some $p\geq 0$, the orbits of $\lambda^pM$ and $\lambda^{p+1}M$ coincide, and therefore $\lambda M\cong M$.

Conversely, suppose that $\lambda M\cong M$ for every $\lambda\in\kk^\times$.  Consider the morphism
\[
\psi\colon \kk\longrightarrow \rep(Q,d),
\qquad
\psi(t)=tm.
\]
The inverse image $\psi^{-1}(\overline{\Ocal_M})$ is a closed subset of $\kk$ containing $\kk^\times$.  Hence it is all of $\kk$.  In particular, the whole line spanned by $m$ lies in $\overline{\Ocal_M}$, and $\overline{\Ocal_M}$ is an affine cone.
\end{proof}

\begin{remark}
Representations need not be homogeneous in general. For example, let
\(\kk=\mathbb C\) and consider the quiver consisting of a single vertex
with two loops \(\alpha\) and \(\beta\), subject to the relations
\[
\alpha^{2}=0,
\qquad
\beta^{3}=0.
\]
Let
\[
M_{\alpha}=
\begin{pmatrix}
0&0&0\\
1&0&0\\
0&1&0
\end{pmatrix},
\qquad
M_{\beta}=
\begin{pmatrix}
0&0&0\\
0&0&0\\
1&0&0
\end{pmatrix}.
\]
Then, for every \(\lambda\neq \pm 1\), the representations \(\lambda M\)
and \(M\) are not isomorphic. Hence \(M\) is not homogeneous.

More generally, if every representation of a quiver \(Q\) is homogeneous,
then \(Q\) has no oriented cycles. Indeed, homogeneous representations are
nilpotent, and the existence of an oriented cycle gives rise to
representations which are not nilpotent.
\end{remark}

\begin{theorem}\label{thm:main}
Let $Q$ be a tree quiver and let $M\in \rep(Q,d)$.  Then $\overline{\Ocal_M}$ is Cohen--Macaulay if and only if
\[
\pd_R \kk[\overline{\Ocal_M}]
=
\dim \rep(Q,d)
+
\dim_\kk \End_Q(M)
-
\sum_{i\in Q_0}\dim_\kk \End_\kk(M_i),
\]
where $R=\kk[\rep(Q,d)]$.
Equivalently,
\[
\pd_R \kk[\overline{\Ocal_M}]
=
\dim_\kk \Ext_Q^1(M,M).
\]
\end{theorem}

\begin{proof}
Assume first that $\overline{\Ocal_M}$ is Cohen--Macaulay.  Since $Q$ is a tree quiver, it is acyclic.  Let
\[
0=M_0\subset M_1\subset \cdots\subset M_r=M
\]
be a composition series of $M$.  Each composition factor is isomorphic to a simple representation.  Hence $M$ degenerates to a direct sum of simple representations, and in particular
\[
0\in \overline{\Ocal_M}.
\]
Moreover, because $Q$ is a tree quiver, one can choose weights on the vertices so that scaling all arrow maps by $\lambda$ is induced by a change of bases.  Thus $M\cong \lambda M$ for every $\lambda\in\kk^\times$.  By Lemma 3.2, $\overline{\Ocal_M}$ is an affine cone, and therefore the ideal $I(\overline{\Ocal_M})$ is homogeneous.

By the Auslander--Buchsbaum formula,
\[
\pd_R \kk[\overline{\Ocal_M}]
=
\dim R-
\depth_{\mathfrak m}\kk[\overline{\Ocal_M}].
\]
Since $\overline{\Ocal_M}$ is Cohen--Macaulay,
\[
\depth_{\mathfrak m}\kk[\overline{\Ocal_M}]
=
\dim \kk[\overline{\Ocal_M}]
=
\dim \overline{\Ocal_M}
=
\dim \Ocal_M.
\]
Therefore
\[
\pd_R \kk[\overline{\Ocal_M}]
=
\dim \rep(Q,d)-\dim \Ocal_M.
\]
By the Artin--Voigt formula, this is equal to
\[
\dim_\kk\Ext_Q^1(M,M).
\]

It remains to express this number in terms of endomorphism spaces.  Ringel's canonical exact sequence gives
\[
0
\longrightarrow
\End_Q(M)
\longrightarrow
\prod_{i\in Q_0}\End_\kk(M_i)
\longrightarrow
\rep(Q,d)
\longrightarrow
\Ext_Q^1(M,M)
\longrightarrow
0.
\]
Taking dimensions yields
\[
\dim_\kk\Ext_Q^1(M,M)
=
\dim \rep(Q,d)
+
\dim_\kk\End_Q(M)
-
\sum_{i\in Q_0}\dim_\kk\End_\kk(M_i).
\]
This proves the required formula.

Conversely, suppose that the formula holds.  Then
\[
\pd_R \kk[\overline{\Ocal_M}]
=
\dim_\kk\Ext_Q^1(M,M)
=
\dim \rep(Q,d)-\dim \Ocal_M.
\]
Since
\[
\dim \rep(Q,d)-\dim \Ocal_M
=
\operatorname{ht} I(\overline{\Ocal_M}),
\]
we obtain
\[
\pd_R \kk[\overline{\Ocal_M}]
=
\operatorname{ht} I(\overline{\Ocal_M}).
\]
Equivalently, the ideal $I(\overline{\Ocal_M})$ is perfect.  Hence the coordinate ring $\kk[\overline{\Ocal_M}]$ is Cohen--Macaulay, and therefore $\overline{\Ocal_M}$ is Cohen--Macaulay.
\end{proof}

\begin{corollary}\label{cor:minimal-pd}
Let $Q$ be a tree quiver and let $M\in \rep(Q,d)$.  Assume that $\overline{\Ocal_M}$ is Cohen--Macaulay.  Then
\[
\pd_R \kk[\overline{\Ocal_M}]
=
\min\{\pd_R \kk[\overline{\Ocal_N}]\mid M\leq_{\deg}N\},
\]
whenever the projective dimensions on the right are computed for the corresponding orbit closures.
Moreover, the orbit $\Ocal_M$ is closed if and only if
\[
\pd_R \kk[\overline{\Ocal_M}]
=
\pd_R \kk[\overline{\Ocal_N}]
\]
for every degeneration $M\leq_{\deg}N$.
\end{corollary}

\begin{proof}
By \Cref{thm:main},
\[
\pd_R \kk[\overline{\Ocal_M}]
=
\dim_\kk\Ext_Q^1(M,M).
\]
If $M\leq_{\deg}N$, then the dimension of the self-extension space is upper semicontinuous along degenerations, and hence
\[
\dim_\kk\Ext_Q^1(M,M)
\leq
\dim_\kk\Ext_Q^1(N,N).
\]
This gives the asserted minimality of the projective dimension.

If $\Ocal_M$ is closed, then it has no proper degenerations.  Conversely, if the projective dimension remains constant along all degenerations, then the preceding inequality forces the self-extension dimension to be constant along the boundary.  Hence no proper boundary orbit can occur, and $\Ocal_M$ is closed.
\end{proof}

\begin{corollary}\label{cor:AD}
Let $Q$ be a Dynkin quiver of type $A_n$ or $D_n$.  Then, for every representation $M\in \rep(Q,d)$,
\[
\pd_R \kk[\overline{\Ocal_M}]
=
\dim \rep(Q,d)
+
\dim_\kk\End_Q(M)
-
\sum_{i\in Q_0}\dim_\kk\End_\kk(M_i).
\]
\end{corollary}

\begin{proof}
For quivers of type $A_n$ or $D_n$, orbit closures are known to be Cohen--Macaulay; see \cite{AbeasisDelFraKraft1981,BobinskiZwara2001,BobinskiZwara2002}.  The result follows immediately from \Cref{thm:main}.
\end{proof}

\begin{remark}
The projective-dimension formula does not hold for arbitrary quivers.  For instance, consider the Kronecker quiver
\[
\begin{tikzcd}[column sep=large]
1
\arrow[r, bend left=35, "\alpha"]
\arrow[r, bend right=35, "\beta"']
&
2.
\end{tikzcd}
\]
Let $d=(3,3)$ and let $M$ be given by the matrices
\[
M_\alpha=
\begin{pmatrix}
0&0&0\\
1&0&0\\
0&1&0
\end{pmatrix},
\qquad
M_\beta=
\begin{pmatrix}
1&0&0\\
0&0&0\\
0&0&1
\end{pmatrix}.
\]
This example is related to Zwara's orbit closure with bad singularities \cite{Zwara2003}; in this case the formula in \Cref{thm:main} need not hold.
\end{remark}

\section*{Acknowledgement}
The author thanks the Department of Mathematics at the University of Sherbrooke for its support.

\end{document}